\documentclass[12pt,reqno]{amsart}

\usepackage{amsmath,amsthm}
\usepackage{tikz}
\usetikzlibrary{calc, arrows.meta, positioning}
\usepackage{amssymb}
\usepackage[english]{babel}
\usepackage{comment}

\usepackage{subcaption}
\usepackage{times}
\usepackage[margin=1.2in]{geometry}
\usepackage{mathrsfs}
\usepackage{hyperref}

\newtheorem{thm}{Theorem}[section]

\newtheorem{lem}[thm]{Lemma}
\newtheorem{prop}[thm]{Proposition}

\newcommand{\pf}{\noindent\textbf{Proof.}$\ $}
\newcommand{\zb}{$\hfill\Box$}
\allowdisplaybreaks

\makeatletter\@addtoreset{equation}{section} \makeatother

\newcommand{\fib}{F}

\begin{document}
	
	\title [Diophantine analysis and the Braid group ${\bf B}_3$]{Diophantine analysis and the Braid group ${\bf B}_3$}

	\author[W. He]{Wei He}
	\address{Wei He (corresponding author): School of Mathematics, Southeast University,
		Nanjing, Jiangsu 211189, China.} \email{hewei@seu.edu.cn}

	\author[W. Lu]{Wenhao Lu}
	\address{Wenhao Lu: School of Mathematics, Southeast University,
		Nanjing, Jiangsu 211189, China.} \email{whlu@seu.edu.cn}
		
	\author[H. Yang]{Hang Yang}
	\address{Hang Yang: School of Mathematics, Southeast University,
		Nanjing, Jiangsu 211189, China.} \email{yh200010@seu.edu.cn}	
	
	\author[R. Yang]{Rongwei Yang}
	\address{Rongwei Yang: Department of Mathematics and Statistics, University at Albany, the State University of New York,
    	Albany, NY 12222, U.S.A.} \email{ryang@albany.edu}

	\maketitle

\begin{abstract}

Given a finite dimensional representation $\pi$ of a finitely generated group $G=\langle g_1, \ldots, g_n\rangle$, the associated characteristic polynomial is defined as $Q_\pi(z):=\det(z_0I+z_1\pi(g_1)+\cdots +z_n\pi(g_n))$, and it is known to contain a good amount of structural information about $G$ and $\pi$. 
This paper is a part of an ongoing project to investigate the number-theoretic properties of the algebraic varieties (called {\em eigensurfaces}) $\{z\in \mathbb{C}^{n+1}: Q_\pi(z)=0\}$. Its focus is the distribution of prime triples in the eigensurface ${\mathcal S}:=\{z\in \mathbb{C}^3: (z_0+z_1+z_2)^2+z_0z_1=0\}$ associated with the braid group ${\bf B}_3$ and its reduced Burau representation. We prove that such triples occur with higher frequency on $\mathcal{S}$ than in the ambient lattice, revealing an unexpected connection between group representation theory and analytic number theory.\\
	
\noindent Key words and phrases: characteristic polynomial, eigensurface, prime triples, asymptotic density, braid group, reduced Burau representation.

\noindent Mathematics Subject Classification (2020): 11D45; 20C30; 47A13.

		
\end{abstract}

\section{Introduction}

Let $\mathbb{Z}$ be the set of integers and $\mathbb{Z}_+$ the set of positive integers. In this paper, the set of prime numbers consists of all positive primes in the ordinary sense together with their negatives. A triple $z=(z_0,z_1,z_2)\in \mathbb{Z}^3$ is said to be \emph{primitive} (or {\em coprime}) if the positive greatest common divisor of its components, denoted by $\gcd z$, equals $1$. A primitive triple $z$ is said to be \emph{prime} if at least one of its coordinates is prime, for example $(-3, 8, 8)$ and $(5, 0, -6)$.

This paper investigates the asymptotic density of prime triples among primitive integer lattice points on algebraic surfaces in $\mathbb{C}^3$. Clearly, most algebraic surfaces contain no positive integer triples at all, for instance the surface $\{z_0^2=\sqrt{2}z_1z_2\}$ or $\{z_0^3+z_1^3=z_2^3\}$, let alone prime triples. Therefore, to make the study meaningful, it is essential to select an appropriate class of surfaces. Here we consider the surface
\begin{align}\label{surfaceequ}
	\mathcal{S}=\{(z_0 + z_1 + z_2)^2 + z_0 z_1 = 0\},
\end{align}
which arises as the {\em projective spectrum} (or {\em eigensurface}) associated with the reduced Burau representation of the braid group ${\bf B}_3$.

\subsection{Background and objective}
The notion of projective spectrum was introduced by the last author in \cite{Ya1}. Let $\mathcal{B}$ be a unital Banach algebra and $\mathbb{A}= (A_0, A_1, \dots, A_n)$ a tuple of elements in $\mathcal{B}$. For $z= (z_0, \dots, z_n) \in \mathbb{C}^{n+1}$, set $\mathbb{A}(z) = z_0A_0 + \dots + z_nA_n$. The {\em projective spectrum} of $\mathbb{A}$ is defined as
\[
P(\mathbb{A})= \{z \in \mathbb{C}^{n+1} \mid \mathbb{A}(z) \text{ is not invertible in } \mathcal{B}\}.
\]
When $\mathcal{B}$ is a matrix algebra, $P(\mathbb{A})$ is simply the algebraic variety $\{\det \mathbb{A}(z)=0\}$. In particular, the polynomial $Q_{\mathbb{A}}(z):=\det (z_0I+z_1A_1+\cdots +z_n A_n)$ is called the ({\em joint}) {\em characteristic polynomial} of the matrices, and its zero variety $Z(Q_{\mathbb{A}})$ is called the corresponding {\em eigensurface}. In the case where $A_1,\dots,A_n$ represent a generating set of a group $G$, the projective spectrum $P(\mathbb{A})$ captures deep intrinsic properties of $G$ and the representation. For instance, for certain self-similar groups, such as the infinite dihedral group, the lamplighter group, and the Grigorchuk group of intermediate growth, the projective spectra have been shown to be closely related to the Julia sets of associated spectral dynamical systems \cite{GY, GrY, ZYL}. In a different direction, characteristic polynomials for finite-dimensional Lie algebras have also been actively studied \cite{AY, GLW, HuY}. For further details, we refer the reader to the monograph \cite{Ya2} and the references therein.

This paper is a part of an ongoing project to investigate the number-theoretic properties of eigensurfaces associated with finite or finitely generated groups. In a previous work \cite{GHY}, we studied the symmetric group ${\bf S}_3$, whose eigensurface under a two-dimensional irreducible representation is $\mathcal{A}=\{z_0^{2} - z_1^{2} + z_2^{2} - z_0z_2=0\}$. It was proved that the density of prime triples among primitive lattice points in $\mathcal{A}\cap\mathbb{Z}_+^3$ exceeds that in $\mathbb{Z}_+^3$, indicating that the eigensurface $\mathcal{A}$ meets primes more frequently. The present paper is the second attempt in this direction, focusing on the braid group ${\bf B}_3= \langle b_1, b_2 \mid b_1b_2b_1 = b_2b_1b_2 \rangle$ and its reduced Burau representation. The situation here turns out to be considerably more intricate.

 The reduced Burau representation $\rho$ \cite{Bu} is the group homomorphism ${\bf B}_3 \to \operatorname{GL}_2(\mathbb{C})$ given by
\[
\rho(b_1) = \begin{pmatrix} -1 & 1 \\ 0 & 1 \end{pmatrix}, \qquad 
\rho(b_2) = \begin{pmatrix} 1 & 0 \\ 1 & -1 \end{pmatrix}.
\]
The characteristic polynomial for $(\rho(b_1), \rho(b_2))$ is
\[
Q(z)=\det (z_0I + z_1\rho(b_1) + z_2\rho(b_2)) = (z_0 + z_1 + z_2)^2 + z_0 z_1,
\]
and the corresponding eigensurface is given in \eqref{surfaceequ}. Numerical experiments confirm that among tens of millions of primitive triples in $\mathcal{S}\cap\mathbb{Z}^3$, the density of prime triples exceeds that in the full lattice $\mathbb{Z}^3$. It is an intriguing question whether this phenomenon holds asymptotically, and we shall conduct an in-depth investigation in this paper.

\subsection{Problem formulation and main results}
Let $\# X$ denote the cardinality of a finite set $X$. For a positive number $N$, set $\pi(N) = \#\{0<p \le N \mid p \text{ is prime}\}$. The Prime Number Theorem \cite[Theorem 6]{HW} states that $\pi(N) = \frac{N}{\log N}(1+o(1))$ as $N\to\infty$, so the density of primes in positive integers less than $N$ is asymptotically $1/\log N$. To formulate a higher-dimensional analogue, we choose a suitable domain in $\mathbb{C}^3$. For $N>0$ and a constant $0<c\le 1$, define
\[
D_c(N)=\{z=(z_0,z_1,z_2)\in\mathbb{Z}^3 \mid |z_0|\le N,\ |z_1|\le N,\ |z_2|\le cN\}.
\]
Although the role of the parameter $c$ appears obscure at this time, subsequent discussion will show that its value does affect the asymptotic density of prime triples in the surface ${\mathcal S}$. Let
\begin{align*}
	X_c(N) &= \{z\in D_c(N) \mid z \text{ is a primitive triple}\},\\
	Y_c(N) &= \{z\in D_c(N) \mid z \text{ is a prime triple}\},
\end{align*}
and set
\[
\mathbb{P}_c(N) = \frac{\# Y_c(N)}{\# X_c(N)}.
\]
Intuitively, the ratio $\mathbb{P}_c(N)$ is the density of prime triples among primitive triples in $D_c(N)$. In Section 2 we estimate $\mathbb{P}_c(N)$ and show that it is asymptotically $3\zeta(3)/\log N$, which is independent of $c$. The situation is quite different for the eigensurface $\mathcal{S}$. Define
\[ X_c^\mathcal{S}(N)= {\mathcal S}\cap X_c(N), \qquad Y_c^\mathcal{S}(N)= \mathcal{S} \cap  Y_c(N),\]
and  set
\[\mathbb{P}_c^\mathcal{S}(N) = \frac{\# Y_c^\mathcal{S}(N)}{\# X_c^\mathcal{S}(N)}.\]
Sections 3 and 4 are devoted to estimating $X_c^{\mathcal{S}}(N)$ and $Y_c^{\mathcal{S}}(N)$, respectively. The estimation of $\#Y_c^{\mathcal{S}}(N)$ is the core of the paper. Our results show that $\mathbb{P}_c^{\mathcal{S}}(N)$ depends on $c$. For the specific value $c_0 = 5/49$, we obtain the following asymptotic bounds, which constitute Theorem \ref{main1}:
\begin{align*}
\liminf_{N\to\infty} \frac{\mathbb{P}_{c_0}^{\mathcal{S}}(N)}{\mathbb{P}_{c_0}(N)} \ge 1.0037,\qquad
\limsup_{N\to\infty} \frac{\mathbb{P}_{c_0}^{\mathcal{S}}(N)}{\mathbb{P}_{c_0}(N)} \le 1.2350.
\end{align*}
	
We do not know whether the limit of $\frac{\mathbb{P}_{c_0}^{\mathcal{S}}(N)}{\mathbb{P}_{c_0}(N)}$ exists as $N\to\infty$. However, somewhat surprisingly, the following limit does exist (Theorem \ref{main2}):
\begin{align*}
	\lim_{c\to 0}\liminf_{N\to \infty} \frac{\mathbb{P}_{c}^\mathcal{S}(N)}{\mathbb{P}_{c}(N)}=\lim_{c\to 0}\limsup_{N\to \infty} \frac{\mathbb{P}_{c}^\mathcal{S}(N)}{\mathbb{P}_{c}(N)}=\frac{\sqrt{5}\,\pi^2}{36\,\zeta(3)\,\log\beta}\approx 1.0598,
\end{align*}
where $\beta = (1+\sqrt{5})/2$.

Thus, prime triples appear more frequently on the eigensurface $\mathcal{S}$ than in the full integer lattice. Together with the results on the symmetric group ${\bf S}_3$ \cite{GHY}, this work reveals a novel connection between group theory, projective spectrum, and analytic number theory. Identifying the key structural ingredients responsible for this enhanced density of primes remains an interesting direction of research.\\

\section{The estimate of $\mathbb{P}_c(N)$}

\par This section gives an estimate of $\mathbb{P}_c(N)$. First, we recall some terminologies and notations which will be used in the sequel. Given two real functions $f$ and $g$ defined on $\mathbb{R}$, with $g$ non-vanishing, we write $f \sim g$ if $\lim_{x\to\infty}f(x)/g(x)= 1$. For $g$ non-negative, the notation $f = O(g)$ means that there exists a constant $C>0$ such that $\lvert f(x) \rvert \leq Cg(x),\ x\in \mathbb{R}$. Moreover, we write $f = o(g)$ in the case $\lim_{x\to\infty}f(x)/g(x)= 0$. The symbols $ \left  \lfloor\cdot \right  \rfloor$ and $\lceil \cdot \rceil$ denote the floor function and the ceiling function, respectively. 

Recall that the M\"obius function $\mu: {\mathbb Z}_+\to \{0, \pm 1\}$ is defined by
\begin{equation}\label{Mobius}
	\mu(d)=\begin{cases}
		1, &\text{if}\ d=1,\\	
		(-1)^k, &\text{if}\ d\ \text{is a product of}\ k\ \text{distinct prime numbers},\\
		0, &\text{if}\  d \ \text{has one or more repeated prime factors}.
	\end{cases}
\end{equation}
Using the M\"obius inversion formula, Nymann \cite{Ny} shows that $k$ positive integers, chosen at random from the set $\{1, \dots, N\}$, are relatively prime with probability $1/\zeta(k)+O(1/N)$ for $k\geq 3$ and probability $1/\zeta(2)+O(\log N/N)$ for $k=2$, where $\zeta$ is the Riemann zeta function. In particular,
\[\zeta(3) = \sum_{n = 1}^{\infty}\frac{1}{n^3} = 1.202056903 \cdots.\]
Nymann's approach enables one to give the following estimate of $\mathbb{P}_c(N)$. It is worth noting that the estimate is independent of $c$.

\begin{thm}\label{generalcase}
	The density
	\[\mathbb{P}_c(N) = \frac{3\zeta(3)}{\log N}(1 + o(1)).\]
\end{thm}
\pf  We estimate $\# X_c(N)$ and $\# Y_c(N)$ separately. Recall that $X_c(N)$ is the set of all primitive triples in $D_c(N)$. To count the set $X_c(N)$, we consider the following subsets of $X_c(N)$. Let 
\begin{align*}
	&X_1(N)=\{z\in X_c(N)\mid 1\le z_0, z_1\le N, 1\le z_2\le cN\},\\
	&X_{21}(N)=\{z\in X_c(N)\mid z_0=0, 1\le z_1\le N, 1\le z_2\le cN\},\\
	&X_{22}(N)=\{z\in X_c(N)\mid 1\le z_0\le N, z_1=0, 1\le z_2\le cN\},\\
	&X_{23}(N)=\{z\in X_c(N)\mid 1\le z_0,z_1\le N, z_2=0\},\\
	&X_3(N)=\{(1,0,0)\}.
\end{align*}
By the symmetry of coordinates, we have
\begin{align*}
	\# X_c(N)=8\# X_1(N)+4\left(\# X_{21}(N)+\# X_{22}(N)+\# X_{23}(N)\right)+ 6\# X_3(N).
\end{align*}
First, it follows from Nymann \cite{Ny} that
\begin{align*}
	\# X_1(N) &= \sum_{d=1}^{N} \mu(d) \cdot  \left  \lfloor \frac{N}{d} \right  \rfloor^2 \cdot  \left  \lfloor\frac{cN}{d} \right  \rfloor
	= \frac{cN^3}{\zeta(3)} + o(N^3),\\
	\# X_{21}(N) &= \sum_{d=1}^{N} \mu(d) \cdot  \left  \lfloor \frac{N}{d} \right  \rfloor \cdot  \left  \lfloor\frac{cN}{d} \right  \rfloor
	= \frac{cN^2}{\zeta(2)} + o(N^2)=o(N^3).\\
\end{align*}
Similarly, we have $\# X_{22}(N)=\# X_{23}(N)=o(N^3)$. Therefore, 
\begin{align*}
	\# X_c(N)=\frac{8cN^3}{\zeta(3)}+o(N^3).
\end{align*}

\par Next, we count the set $Y_c(N)$, which is the set of prime triples in $X_c(N)$. Let
\begin{align*}
	Y_*(N)=\{z\in X_*(N) \mid z\ \text{is a prime triple}\},
\end{align*}                                                         where $*$ stands for the subscript of the corresponding set.        
The set $Y_1(N)$ consists of all triples of positive integers in $D_c(N)$ which are primitive and contain at least one prime. By the Prime Number Theorem, the number of positive integer triples in $D_c(N)$ with at least one prime component is 
\[N^2\lfloor cN \rfloor-(N-\pi(N))^2 (\lfloor cN \rfloor-\pi(cN))=\frac{3cN^3}{\log N} + o\left(\frac{N^3}{\log N}\right).\]
We then subtract the number of such triples which are not primitive. But this occurs only if one of coordinates $z_0, z_1, z_2$ is a prime $p$ and the other two are integer multiples of $p$. Therefore, this number is
\begin{align*}
	&2\sum_{\substack{1\le p\le N\\ p\ \text{is prime}}}\left( \left  \lfloor\frac{N}{p} \right  \rfloor  \left  \lfloor\frac{cN}{p} \right  \rfloor \right) 
	+\sum_{\substack{1\le p\le cN\\ p\ \text{is prime}}}\left( \left  \lfloor\frac{N}{p} \right  \rfloor \right)^2 \\
	\leq{} &3\sum_{p \leq N }\left(\frac{N}{p}\right)^2 \leq O(N^2) = o\left(\frac{N^3}{\log N}\right).
\end{align*}
\noindent Hence 
\[\# Y_1(N) = \frac{3cN^3}{\log N} + o\left(\frac{N^3}{\log N}\right). \]	
Observe that for $1\leq j\leq 3$,
\[\# Y_{2j}(N) \le X_{2j}(N) = o\left(\frac{N^3}{\log N}\right), \]	
and $\# Y_3(N)=0$. It follows that 
\begin{align*}
	\# Y_c(N)&=8\# Y_1(N)+4\left(\# Y_{21}(N)+\# Y_{22}(N)+\# Y_{23}(N)\right)\\
	&=\frac{24cN^3}{\log N} + o\left(\frac{N^3}{\log N}\right).
\end{align*}

In view of the estimates of $\# X_c(N)$ and $\# Y_c(N)$ above, we have \[\mathbb{P}_c(N)= \frac{\# Y_c(N)}{\# X_c(N)} = \frac{3\zeta(3)}{\log N}(1 + o(1)),\]
completing the proof of the theorem. \zb\\

\section{Parameterization of $\mathcal{S}$ and an estimate of $\# X_c^\mathcal{S}(N)$}

This section gives an estimate of $\#X_c^\mathcal{S}(N)$, where 
\begin{align*}  
	X_c^\mathcal{S}(N)= \{z\in \mathcal{S} \mid \gcd{z}= 1\} \cap D_c(N).
\end{align*}
To this end, we first parameterize coprime triples in $\mathcal{S}$. Recall that $\mathcal{S}=\{(z_0 + z_1 + z_2)^2 + z_0 z_1 = 0\}$. A direct verification shows that if $z$ is a coprime triple in $\mathcal{S}$,  then all the coordinates of $z$ are non-zero, except for the four triples $(0,1,-1)$, $(0,-1,1)$, $(1,0,-1)$, $(-1,0,1)$. Hence in the sequel, we only count the number of primitive triples in $\mathcal{S}$ with all coordinates non-zero. Let  
\begin{align*}  
	&\Delta^\mathcal{S}=\{z\in \mathcal{S} \mid \gcd{z}= 1\}, \\
	&\phi:\mathbb{Z}^2\to \mathbb{Z}^3, (m, n)\mapsto (-m^2, n^2, m^2+mn-n^2).
\end{align*}
The map $\phi$ has the following property.
\begin{lem}\label{mncoprime}
	For $m,n\in \mathbb{Z}$, we have $\gcd(m, n)=1$ if and only if $\gcd \phi(m,n)=1$.
\end{lem}

\pf If $\gcd(m, n)=1$, then by the fundamental theorem of arithmetic, we have $\gcd (-m^2, n^2)=1$, and thus $\gcd \phi(m,n)=\gcd (-m^2, n^2, m^2+mn-n^2)=1$.

Conversely, if $\gcd(m, n) \neq 1$, then $\gcd(m, n)=d>1$. Since $d \mid m$, $d \mid n$, we have $d\mid \gcd \phi(m,n)$, 
completing the proof.\zb\\

We define $\Delta_{\pm}=\{\pm \phi(m,n) \mid \gcd(m,n)=1\}$. The next lemma gives a parameterization of the primitive triples in $\mathcal{S}$.

\begin{lem}\label{parameterize} $\Delta^\mathcal{S}=\Delta_+\cup \Delta_-$.\end{lem}

\pf Since $\mathcal{S}=\{(z_0 + z_1 + z_2)^2 + z_0 z_1 = 0\}$, it is easy to verify that $\pm \phi(m,n)\in \mathcal{S}$ for any integers $m,n$. In view of Lemma \ref{mncoprime}, we have $\Delta_+\cup \Delta_-\subseteq \Delta^\mathcal{S}$. 
So we only need to show that $\Delta^\mathcal{S}\subseteq \Delta_+\cup \Delta_-$.

Let $z=(z_0, z_1, z_2)\in \Delta^\mathcal{S}$, then $(z_0 + z_1 + z_2)^2 + z_0 z_1 = 0$ and $\gcd (z_0, z_1, z_2)=1$. We first show that $\gcd (z_0, z_1)=1$. In fact, if there exists a prime $p$ such that $p\mid z_0$, $p\mid z_1$, then $p^2\mid (z_0 + z_1 + z_2)^2$. It follows that $p\mid z_0 + z_1 + z_2$, and thus $p\mid z_2$, contradicting with the fact that  $\gcd (z_0, z_1, z_2)=1$.

By the fact that $\gcd (z_0, z_1)=1$ and $(z_0 + z_1 + z_2)^2=-z_0 z_1$, we know that both $|z_0|$ and $|z_1|$ are perfect squares, and $z_0 z_1 < 0$. If $z_0 < 0$, then $z_1 > 0$. In this case, there exist integers $m,n$ such that $z_0=-m^2$, $z_1=n^2$, and $z_2=m^2+mn-n^2$. Thus $z\in \Delta_+$. If $z_0 > 0$, then $z_1 < 0$. In this case, there exist integers $m',n'$ such that $z_0=m'^2$, $z_1=-n'^2$, and $z_2=-m'^2-m'n'+n'^2$. Thus $z\in \Delta_-$. Hence $\Delta^\mathcal{S}\subseteq \Delta_+\cup \Delta_-$, completing the proof.\zb\\                                                       

To simplify the computation, we make use of the symmetries of sets. Let 
\begin{align*}
	&E_c(N)=\{(m,n)\in \mathbb{Z}^2 \mid 0< m\le N, 0<|n|\le N , 0< m^2+mn-n^2\le cN^2\},\\
	&E'_c(N)=\{(m,n)\in \mathbb{Z}^2 \mid 0< m\le N, 0<|n|\le N, -cN^2\le m^2+mn-n^2< 0\}.
\end{align*}
Let $\operatorname{Sgn}(n)$ denote the sign of the integer $n$. Then the map $(m,n)\mapsto \operatorname{Sgn}(n)(n,-m)$
establishes a one-to-one correspondence between $E_c(N)$ and $E'_c(N)$, and hence $\# E_c(N)=\# E'_c(N)$. The diagram for $E_c(N)$ is shown in Figure \ref{fig:1}.

Furthermore, we consider the following two subsets of $E_c(N)$:
\begin{align*}
	&\Omega_c^X(N):=\{(m,n) \in E_c(N)\mid \gcd(m,n)=1\},\\
	&\Omega_c^Y(N):=\{(m,n) \in E_c(N)\mid m^2+mn-n^2\ \text{is a prime}\}.
\end{align*}
\begin{lem}\label{countmn}	
$\# X_c^\mathcal{S}(N^2)=4\# \Omega_c^X(N)$, $\# Y_c^\mathcal{S}(N^2)=4\# \Omega_c^Y(N)$.
\end{lem}
\pf Clearly, points of $\Delta_+$ and $\Delta_-$ are in one-to-one correspondence with each other. By Lemma \ref{parameterize},
\begin{align*}
 \# X_c^\mathcal{S}(N^2)=\# (\Delta^\mathcal{S}\cap D_c(N^2))=2\# (\Delta_+\cap D_c(N^2)).
\end{align*}
Observe that the map $\phi$ restricted to $\mathbb{Z}_+\times\mathbb{Z} $ is injective, and hence
\begin{align*}
\# (\Delta_+\cap D_c(N^2))&=\#\{\phi(m,n)\in D_c(N^2)|~\gcd(m,n)=1\} \\
&=\# \{(m,n)\in E_c(N)\cup E'_c(N)|~\gcd(m,n)=1\}\\
&=2\#\Omega_c^X(N).
\end{align*}
It follows that $\# X_c^\mathcal{S}(N^2)=4\# \Omega_c^X(N)$.

To count the set $Y_c^\mathcal{S}(N^2)$, observe that $\phi(m,n)$ is a prime triple if and only if $\gcd(m,n)=1$ and $m^2+mn-n^2$ is a prime. Similar reasoning as above gives that $\# Y_c^\mathcal{S}(N^2)=4\# \Omega_c^Y(N)$. \zb\\

The parameterization $\phi$ of $\mathcal{S}$ and the preceding discussions lead us to the following expression of $\mathbb{P}_c^\mathcal{S}(N^2)$.

\begin{thm}\label{PN2}
	It holds that 
	\begin{align*}
		\mathbb{P}_c^\mathcal{S}(N^2) =\frac{\# \Omega_c^Y(N)}{\# \Omega_c^X(N)}.
	\end{align*}
\end{thm}

To compute $\mathbb{P}_c^\mathcal{S}(N^2)$, we will first estimate $\# \Omega_c^X(N)$. The estimate of $\# \Omega_c^Y(N)$ is left to the next section.

\begin{prop}\label{XN}
The cardinality of $\Omega_c^X(N)$ is given by
\begin{align*}
		\# \Omega_c^X(N)=\begin{cases}			
		\frac{6}{\pi^2}\left(\frac{1}{2} +\frac{1}{\sqrt{5}} \log \frac{3+\sqrt{5}}{2}\right) N^2+o(N^2), &\ \text{if}\ c=1,\\
		\frac{6}{\pi^2}\left(1-\int^1_{\frac{1-\sqrt{5-4c}}{2}}
\left(1-\frac{\sqrt{5n^{2}+4c}-n}{2}\right)dn\right)N^2+o(N^2), &\ \text{if}\ 0<c<1.
		\end{cases}
\end{align*}
\end{prop}

\pf Since $E_c(N)$ is a bounded region in $\mathbb{R}^2$, by \cite{Ro} and \cite[Theorem 459]{HW}, we have that the number of primitive pairs in $E_c(N)$ is approximately $1/\zeta(2)=6/\pi^2$ times the area of the region. Let $\mathcal{A}(c)$ denote the area of $E_c(N)$. From Figure \ref{fig:1}, one sees that 
\begin{align*}
	\mathcal{A}(c)&=N^2-\left(\int^1_{\frac{1-\sqrt{5-4c}}{2}}
	\left(1- \frac{\sqrt{5n^2+4c}-n}{2}\right) dn\right)N^2\\
	&=\begin{cases}			
		\left(\frac{1}{2} +\frac{1}{\sqrt{5}} \log \frac{3+\sqrt{5}}{2}\right) N^2, &\ \text{if}\ c=1,\\
			\left(1-\int^1_{\frac{1-\sqrt{5-4c}}{2}}
		\left(1-\frac{\sqrt{5n^{2}+4c}-n}{2}\right)dn\right)N^2+o(N^2), &\ \text{if}\ 0<c<1.
	\end{cases}	
\end{align*}
and the proposition follows. \zb\\

\section{An estimate of $\# Y_c^\mathcal{S}(N)$}

In this section, we will give an estimate of $\# Y_c^\mathcal{S}(N)$. Recall that 
\begin{align*}  
	Y_c^\mathcal{S}(N)&= \{z\in \mathcal{S} \cap D_c(N) \ |~z\ \text{is a prime triple}\}.
\end{align*}
By Lemma \ref{countmn}, counting $\# Y_c^\mathcal{S}(N^2)$ amounts to computing $\# \Omega_c^Y(N)$, where 
\[\Omega_c^Y(N)=\{(m,n) \in E_c(N)\mid m^2+mn-n^2\ \text{is a prime}\}.\]
Hence we need to count the number of coprime pairs $(m,n)$ such that $0<m\leq N$, $0<|n|\leq N$ and $m^2+mn-n^2$ is a positive prime $p$ not exceeding $cN^2$. Thus the main task of this section is to count the solutions of the Diophantine equation
\begin{equation} \label{eq:main}
	m^2 + mn - n^2 = p.
\end{equation}

\subsection{The structure of solutions to Equation \eqref{eq:main}}

To count $\# \Omega_c^Y(N)$, it suffices to consider the integer solutions $(m,n)$ with $m>0$, and this is what we will do in the sequel. Note that if a pair of integers $(m,n)$ is a solution for Equation \eqref{eq:main}, then $\gcd(m,n)=1$. The converse is not true, for example $(9,7)$. Hence we need to know for which prime numbers $p$, Equation \eqref{eq:main} has integer solutions. The following proposition from \cite[Equation 3.18f]{Co} addresses this issue.

\begin{prop}\label{existence}
For a positive prime $p$, Equation \eqref{eq:main} has integer solutions if and only if \(p = 5\) or \(p \equiv \pm 1 \pmod 5\).
\end{prop}

 Additionally, Cohn \cite[Section 14]{Co} explains the connection between the solvability of \eqref{eq:main} and the principal prime ideals in the ring of algebraic integers $\mathcal{O}_{\mathbb{Q}(\sqrt{5})}$. The Diophantine equation \eqref{eq:main} has also been studied recently in \cite{HS} in a more general setting using a geometric approach. In particular, it shows that, when it has integer solutions, the set of all such solutions has a recursive structure.

\begin{figure}[htbp]
	\centering
	\def\N{3} 
	
	\begin{tikzpicture}[scale=1.2, >=Stealth]
		\def\phiA{1.618}
		\def\phiB{-0.618}
		
		\def\A{0.63245*\N} 
		\def\SqrtFive{2.236068}
		
		\def\tMax{2.0} 
		
		\begin{scope}
			
			\clip (0,-0.8*\N) rectangle (\N,\N);
			
			\fill[gray!30] (0,0) -- (\N, \N*\phiA) -- (\N, \N*\phiB) -- cycle;
			
			\fill[white] 
			plot[variable=\t, domain=\tMax:-\tMax, samples=100] 
			({ \A*cosh(\t) }, { 0.5*(\A*cosh(\t) + \SqrtFive*\A*sinh(\t)) })
			-- cycle; 
			
			\draw[thick, black] (0,0) -- (\N, \N*\phiA);
			\draw[thick, black] (0,0) -- (\N, \N*\phiB);
			
			\draw[thick, black] 
			plot[variable=\t, domain=-\tMax:\tMax, samples=100] 
			({ \A*cosh(\t) }, { 0.5*(\A*cosh(\t) + \SqrtFive*\A*sinh(\t)) });
		\end{scope}
		
		\draw[densely dashed, black, thin] (\N, -0.8*\N) -- (\N, 1.2*\N) 
		node[above, text=black, font=\footnotesize] {$m=N$};
		
		\draw[densely dashed, black, thin] (-0.2*\N, \N) -- (1.2*\N, \N) 
		node[right, text=black, font=\footnotesize] {$n=N$};
		
		\def\mZero{1.4} 
		\def\nZero{0.4}
		\coordinate (P1) at (\mZero, \nZero);
		\coordinate (P2) at (\mZero, \mZero-\nZero);
		\coordinate (MX) at (\mZero, 0);
		
		\draw[dashed, black] (MX) -- (P1);
		
		\fill[black] (P1) circle (1.5pt);
		\node[right, font=\scriptsize, text=black] at (P1) {$(m_0, n_0)$};
		
		\fill[black] (P2) circle (1.5pt);
		\node[left, font=\scriptsize, text=black, yshift=-2pt] at (P2) {$(m_0, m_0-n_0)$};
		
		\node[below, font=\scriptsize] at (MX) {$m_0$};
		
		\draw[densely dotted, darkgray] (0,0) -- (1.2*\N, 1.2*\N) node[right, font=\tiny] {$m=n$};
		\draw[densely dotted, darkgray] (0,0) -- (1.5*\N, 0.75*\N) node[right, font=\tiny] {$m=2n$};
		
		\draw[->, thick] (-0.3*\N,0) -- (1.5*\N,0) node[right] {$m$};
		\draw[->, thick] (0,-0.8*\N) -- (0,1.5*\N) node[above] {$n$};
		
		\node[below left] at (0,0) {$O$};
		\node[below, xshift=-5pt] at (1.15*\N, 0) {$N$}; 
		\node[left, yshift=-5pt] at (0, 1.15*\N) {$N$};
		
		\node[above, font=\tiny, rotate=58] at (0.6*\N, \N) {$m^2+mn-n^2=0$};
		\node[right, font=\tiny] at (0.65*\N, 0.4*\N) {$m^2+mn-n^2=cN^2$};		
	\end{tikzpicture}    
	\caption{The region $E_c(N)$ and location of representative solutions}
	\label{fig:1}
\end{figure}
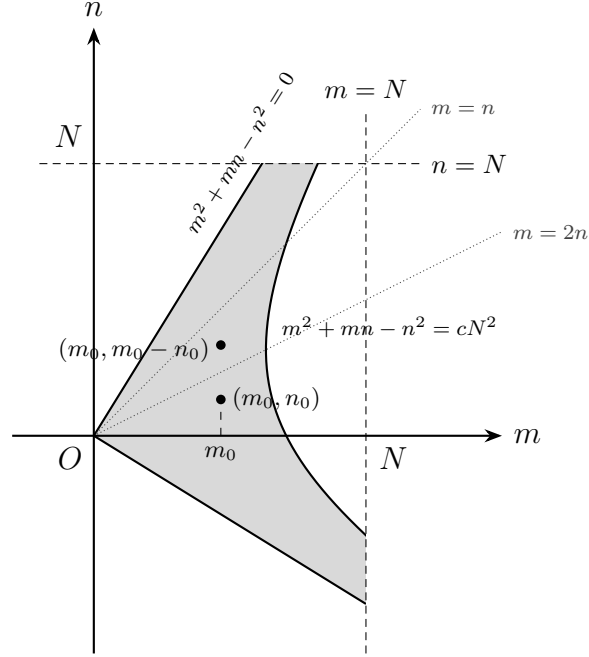

An integer solution $(m_0, n_0)$ of Equation \eqref{eq:main} is said to be \emph{representative}, if  $m_0$ is the smallest among all the integer solutions $(m, n)$ with $m>0$. It is proved \cite[Proposition 4.3]{HS} that, if $p=5$, then $(2,1)$ is the only representative solution for Equation \eqref{eq:main}. While in the case $p \equiv \pm 1 \pmod 5$, there are exactly two 
representative solutions $(m_0, n_0)$ and $(m_0, m_0-n_0)$, satisfying 
\begin{equation}\label{repsolu1}
\begin{aligned}
&(m_0, n_0)\in \left[\sqrt{\frac{4p}{5}}, \sqrt{p}\right]\times \left[0, \sqrt{\frac{p}{5}}\right],\\
&(m_0, m_0-n_0)\in \left[\sqrt{\frac{4p}{5}}, \sqrt{p}\right]\times \left[\sqrt{\frac{p}{5}}, \sqrt{p}\right].
\end{aligned}
\end{equation}
One can also verify that $(m_0, n_0)$ lies between the positive $m$-axis and the line $\{m = 2n\}$, and $(m_0, m_0-n_0)$ lies between the lines $\{m = 2n\}$ and $\{m = n\}$, as shown in Figure \ref{fig:1}.

The general integer solutions of Equation \eqref{eq:main} are generated by the following matrix in the special linear group $SL(2,{\mathbb Z})$: \[A = \begin{pmatrix} 1 & 1 \\ 1 & 2 \end{pmatrix}.\] 

\begin{prop}[\cite{HS}]\label{generalsolu}
A primitive pair $(m,n)$ is an integer solution of Equation \eqref{eq:main} if and only if it is of the form
\begin{align*}
	\begin{pmatrix} m \\ n \end{pmatrix} 
	= A^k \begin{pmatrix} m_0 \\ n_0 \end{pmatrix},\ \ \ k\in {\mathbb Z},
\end{align*}
where $(m_0, n_0)$ is representative.
\end{prop}

 The function
\[\pi(N; q, a) := \#\{0<p \leq  N \ | \ p\ \text{is prime},\ p\equiv a \pmod q\}\]
counts the number of primes $p \equiv a \pmod q$ up to $N$. Landau's theorem gives an estimate of $\pi(N; q, a)$ (\cite[Section 120]{La} and \cite[Theorem 7.25]{LeV}). 

\begin{lem}\label{Lan}
	Let $a, q > 0$ be relatively prime integers. Then 
	\[\pi(N; q, a) = \frac{1}{\varphi(q)} \frac{N}{\log N} \big(1 + o(1)\big),\]
	where $\varphi(q)$ is Euler's totient function.
\end{lem}

For instance, since $\varphi(5)=4$, asymptotically we have \[\pi(N; 5, 1)=\pi(N; 5, 4)=\frac{N}{4\log N}(1+o(1)).\] 

\vskip 3mm

\subsection{The number of solutions to Equation \eqref{eq:main}}\label{sec:solution_count} Proposition \ref{existence} indicates that Equation \eqref{eq:main} has solutions if and only if $p = 5$ or $p \equiv \pm 1 \pmod 5$. But for different prime $p$, the number of solutions within $E_c(N)$ differs. In this subsection, we will give a sharp estimate dependent on $p$. It suffices to do so in $E_1(N)$, since  $E_c(N)\subset E_1(N)$ for $0<c<1$.

For any given prime $p\in[1, N^2]$ such that \eqref{eq:main} has solutions, let \((m_0, n_0)\) be the representative solution to \eqref{eq:main} such that $m_0>2n_0>0$. Then 
\begin{equation} \label{m0n0}
	m_0^2 + m_0 n_0 - n_0^2 = p.
\end{equation}
By Proposition \ref{generalsolu}, all its derived solutions are given by
\begin{align*}
	\begin{pmatrix} m_k \\ n_k \end{pmatrix} = 
	\begin{pmatrix} 1 & 1 \\ 1 & 2 \end{pmatrix}^k 
	\begin{pmatrix} m_0 \\ n_0 \end{pmatrix}, \quad k \in \mathbb{Z}.	
\end{align*}
Let $\beta = (1+\sqrt{5})/2$, $\gamma =-\beta^{-1}=(1-\sqrt{5})/2$. Then 
\begin{align}\label{mknk}
	\begin{pmatrix} m_k \\ n_k \end{pmatrix} = 
	\begin{pmatrix} (\beta^{2k-1}-\gamma^{2k-1})/\sqrt{5} & (\beta^{2k}-\gamma^{2k})/\sqrt{5} \\ (\beta^{2k}-\gamma^{2k})/\sqrt{5} & (\beta^{2k+1}-\gamma^{2k+1})/\sqrt{5} \end{pmatrix} 
	\begin{pmatrix} m_0 \\ n_0 \end{pmatrix}.
\end{align}
When $k>0$, since \(\fib_k=(\beta^k-\gamma^k)/\sqrt{5}\) is the \(k\)-th Fibonacci number, the above formula can be rewritten as 
\begin{align*}
\begin{cases}
	m_k = \fib_{2k-1} m_0 + \fib_{2k} n_0,\\ 
	n_k = \fib_{2k} m_0 + \fib_{2k+1} n_0.
\end{cases}
\end{align*}
It follows that $n_k>m_k>0$ and $n_k$ is monotonically increasing. Let $k_1^+$ be the maximal $k$ such that $(m_k, n_k)\in E_1(N)$ when $0<k\leq k_1^+$. The above reasoning implies that $k_1^+$ is in fact the maximal $k$ satisfying $n_k\leq N$. We solve from Equation \eqref{m0n0} and \eqref{mknk} that 
\begin{equation} \label{k1+}
	k_1^+ = \left\lfloor \frac{\log\left( \sqrt{5}N + \sqrt{5N^2 + 4p} \right) - \log(2\beta n_0 + 2m_0)}{2 \log \beta} \right\rfloor.
\end{equation}
When \(k < 0\), let $k_1^-$ be the maximal absolute value of $k$ such that $(m_k, n_k)\in E_1(N)$ when $-k_1^-<k<0$. Similar argument as above gives that 
\begin{equation} \label{k1-}
	k_1^- =  \left\lfloor \frac{\log\left( \sqrt{5}N + \sqrt{5N^2 - 4p} \right) - \log(2\beta m_0 - 2n_0)}{2 \log \beta} \right\rfloor.
\end{equation} 
For another representative solution $(m_0, m_0 - n_0)$, let 
\begin{align} 
	k_2^+ &= \left\lfloor \frac{\log\left( \sqrt{5}N + \sqrt{5N^2 + 4p} \right) - \log(2(\beta+1)m_0 - 2\beta n_0)}{2 \log \beta} \right\rfloor,\label{k2+}\\	
	k_2^- &= \left\lfloor \frac{\log\left( \sqrt{5}N + \sqrt{5N^2 - 4p} \right) - \log(2(\beta-1)m_0 + 2n_0)}{2 \log \beta} \right\rfloor. \label{k2-}
\end{align}
One can infer that the derived solutions of $(m_0, m_0 - n_0)$ are in $E_1(N)$ if and only if $-k_2^-\leq k\leq k_2^+$. 
Therefore, for fixed prime $p\in [1, N^2]$, $p=\pm 1 \pmod 5$, the total number of solutions within $E_1(N)$ is 
\begin{align*}
	s_p:=k_1^+ + k_1^- + k_2^+ + k_2^- + 2,  
\end{align*}
and for $p=5$, $s_5=k_1^+ + k_1^-+1$.

Now we estimate \(s_p\) for $p$ in different intervals. To this end, we first rewrite $k_1^+$, $k_1^-$, $k_2^+$ and $k_2^-$ in an appropriate form. Let $t_p=(\beta m_0-n_0)/\sqrt{p}$. Since $m_0>2n_0>0$, we have $t_p\in (1, \beta)$. Then it follows from Equation \eqref{m0n0} that 
\begin{align*}
	&2\beta m_0 - 2n_0=2t_p\sqrt{p},\\
	&2\beta n_0 + 2m_0=\frac{(2\beta n_0 + 2m_0)(2\beta m_0 - 2n_0)}{2\beta m_0 - 2n_0}=\frac{2\beta}{t_p}\sqrt{p},\\	
	&2(\beta-1)m_0 + 2n_0=\frac{(2(\beta-1)m_0 + 2n_0)(2\beta m_0 - 2n_0)}{2\beta m_0 - 2n_0}=\frac{2}{t_p}\sqrt{p},\\
	&2(\beta+1)m_0 - 2\beta n_0=\frac{(2(\beta+1)m_0 - 2\beta n_0)(2(\beta-1)m_0 + 2n_0)}{2(\beta-1)m_0 + 2n_0}=2\beta t_p\sqrt{p}.
\end{align*}
Substituting the above formulas into \eqref{k1+}--\eqref{k2-} yield
\begin{equation} \label{ks}
	\begin{aligned}
		k_1^+ &= \left\lfloor \frac{\log\left( \sqrt{\frac{5N^2}{p}} + \sqrt{\frac{5N^2}{p} + 4} \right) - \log\left(\frac{2\beta}{t_p}\right)}{2 \log \beta} \right\rfloor, \\	
		k_1^- &= \left\lfloor \frac{\log\left( \sqrt{\frac{5N^2}{p}} + \sqrt{\frac{5N^2}{p} - 4} \right) - \log(2t_p)}{2 \log \beta} \right\rfloor, \\
		k_2^+ &= \left\lfloor \frac{\log\left( \sqrt{\frac{5N^2}{p}} + \sqrt{\frac{5N^2}{p} + 4} \right) - \log(2\beta t_p)}{2 \log \beta} \right\rfloor, \\
		k_2^- &= \left\lfloor \frac{\log\left( \sqrt{\frac{5N^2}{p}} + \sqrt{\frac{5N^2}{p} - 4} \right) - \log\left(\frac{2}{t_p}\right)}{2 \log \beta} \right\rfloor.
	\end{aligned}
\end{equation}

It is easy to see that $s_5 = O(\log N)$. For $s_p$ where $p=\pm 1 \pmod 5$, we partition the interval $(0, N^2]$ in the following way. Define
\begin{align*}
	F(m, u) &= \frac{5N^2}{(u \beta^{2m} + u^{-1} \beta^{-2m})^2}, \\
	G(m, u) &= \frac{5N^2}{\left(u^{-1}\beta^{2m+1} - u\beta^{-2m-1}\right)^2},
\end{align*}
where \(u \in [\beta^{-1}, \beta]\), and $m$ is in the set of non-negative integers. Since $F(m, \beta)$ is a monotonically decreasing function of $m$, and $F(0, \beta)=N^2$, $F(+\infty, \beta)=0$, we get a partition of 
\((0, N^2]\) by the tagged points $F(m, \beta)$, $m\geq 0$. Let
\begin{align*}
	A_m = (F(m+1, \beta), F(m, \beta)].
\end{align*}
Then the intervals $A_m$ are pairwise disjoint and $(0, N^2] = \bigcup_{m=0}^{\infty} A_m$. 

We further refine the partition by adding three tagged points into each \(A_m\). Let 
\begin{align*}
	u_{m1} = \sqrt{\frac{\beta^{4m+3} - 1}{\beta^{4m+3} + 1}\cdot \frac{1}{\beta}},\quad
	u_{m2} = 1, \quad
	u_{m3}= \sqrt{\frac{\beta^{4m+5} - 1}{\beta^{4m+5} + 1}\cdot \beta}.
\end{align*}
It is worth noting that $u_{m1}$ and $u_{m3}$ are chosen so that 
\begin{align}\label{FG}
	F(m+1, u_{m1}) = G(m, u_{m1}),\quad F(m+1, u_{m3}) = G(m+1, u_{m3}),
\end{align}
One verifies that 
\begin{align}
	\beta^{-1} \leq u_{m1} < u_{m2} < u_{m3} <\beta,
\end{align}
and the equality holds if and only if $m=0$. The equality $\beta^{-1}=u_{01}$ implies that $A_0$ has two tagged points added, not three. But this does not lead to any other difference in the computing process, so we ‌write them down in a uniform way. Since $F(m+1, u)$ is a monotonically decreasing function of $u$ on the interval \([\beta^{-1}, \beta]\), we have 
\begin{align*}
	F(m+1, \beta) < F(m+1, u_{m3}) < F(m+1, u_{m2}) < F(m+1, u_{m1}) \leq F(m+1, \beta^{-1}).
\end{align*}
Observing that $F(m+1, \beta^{-1})= F(m, \beta)$, we get a partition of $A_m$ into four pairwise disjoint subsets
\begin{align*}
	&A_{m1} = \big(F(m+1, \beta),\; F(m+1, u_{m3})\big], \quad 
	A_{m2} = \big(F(m+1, u_{m3}),\; F(m+1, u_{m2})\big], \\	 
	&A_{m3} = \big(F(m+1, u_{m2}),\; F(m+1, u_{m1})\big], \quad
	A_{m4} = \big(F(m+1, u_{m1}),\; F(m, \beta)\big].
\end{align*}
Then we obtain the following estimates of $s_p$ for $p$ in each $A_{mi}$.

\begin{thm}\label{numofsolutions}
For a prime $p$ with $p=\pm 1 \pmod 5$, and the intervals $A_{mi}$ defined above, the following estimates for $s_p$ hold. \\
\emph{(i)} For \(p \in A_{m1}\), \( 4m+4 \leq s_p \leq 4m + 6\);\\
\emph{(ii)} For \(p \in A_{m2}\), \( 4m+3 \leq s_p \leq 4m + 5\);\\
\emph{(iii)} For \(p \in A_{m3}\), \( 4m+2 \leq s_p \leq 4m + 4\);\\
\emph{(iv)} For \(p \in A_{m4}\), \( 4m+1 \leq s_p \leq 4m + 3\).\\
\end{thm}

\pf We only prove (i) in detail. To estimate $s_p$ for \(p \in A_{m1}\), we will estimate $k_1^+$, $k_1^-$, $k_2^+$ and $k_2^-$ respectively. Since \(p \in A_{m1}\), we have 
\begin{align*}
	G(m+1,\beta^{-1})<F(m+1, \beta)<p\leq F(m+1,u_{m3})=G(m+1,u_{m3}),
\end{align*}
where the first inequality is directly verified and the last equality follows from \eqref{FG}. We first estimate $k_1^+$. On the one hand, $p\leq G(m+1,u_{m3})$ and Formula \eqref{ks} for $k_1^+$ imply that 
\begin{align*}
	k_1^+ &= \left\lfloor \frac{\log\left( \sqrt{\frac{5N^2}{p}} + \sqrt{\frac{5N^2}{p} + 4} \right) - \log\left(\frac{2\beta}{t_p}\right)}{2 \log \beta} \right\rfloor\\	
	& \geq \left\lfloor \frac{\log\left( \sqrt{\frac{5N^2}{G(m+1, u_{m3})}} + \sqrt{\frac{5N^2}{G(m+1, u_{m3})} + 4} \right) - \log\left(\frac{2\beta}{t_p}\right)}{2 \log \beta} \right\rfloor\\
	&=\left\lfloor\frac{\log\left(2\beta^{2m+3}/u_{m3}\right) - \log(2\beta / t_p)}{2\log \beta}\right\rfloor
	=\left\lfloor m+1+\frac{\log (t_p/u_{m3})}{2\log \beta}       \right\rfloor.
\end{align*}
Note that $t_p \in (1, \beta)$ and $u_{m3}\in (1, \beta)$, so we have 
\begin{align*}
	k_1^+ \geq 
	\begin{cases}
		m, &\textrm{if}\ t_p \in (1, u_{m3}),\\
		m+1, &\textrm{if}\ t_p \in [u_{m3}, \beta).
	\end{cases}
\end{align*}
On the other hand, the inequality $p>G(m+1, \beta^{-1})$ implies that 
\begin{align*}
	k_1^+	& \leq \left\lfloor \frac{\log\left( \sqrt{\frac{5N^2}{G(m+1, \beta^{-1})}} + \sqrt{\frac{5N^2}{G(m+1, \beta^{-1})} + 4} \right) - \log\left(\frac{2\beta}{t_p}\right)}{2 \log \beta} \right\rfloor\\
	&=\left\lfloor m+1+\frac{\log (t_p\beta)}{2\log \beta}       \right\rfloor\leq m+1.
\end{align*}
To estimate $k_1^{-}$, we substitute the inequality $p\leq F(m+1, u_{m3})$ to Formula \eqref{ks} for $k_1^{-}$. Thus
\begin{align*}
	k_1^- & \geq \left\lfloor \frac{\log\left( \sqrt{\frac{5N^2}{F(m+1, u_{m3})}} + \sqrt{\frac{5N^2}{F(m+1, u_{m3})} - 4} \right) - \log(2t_p)}{2 \log \beta} \right\rfloor\\	
	& = \left\lfloor\frac{\log\left(2u_{m3} \cdot \beta^{2(m+1)}\right) - \log 2t_p}{2\log \beta}\right\rfloor	=\left\lfloor m+1+\frac{\log (u_{m3}/t_p)}{2\log \beta}       \right\rfloor,
\end{align*}
and then 
\begin{align*}
	k_1^- \geq 
	\begin{cases}
		m+1, &\textrm{if}\ t_p \in (1, u_{m3}],\\
		m, &\textrm{if}\ t_p \in (u_{m3}, \beta).
	\end{cases}		
\end{align*}
The upper bound of \( k_1^{-} \) follows from the fact that $p>F(m+1, \beta)$, which implies
\begin{align*}
	k_1^- & \leq \left\lfloor \frac{\log\left( \sqrt{\frac{5N^2}{F(m+1, \beta)}} + \sqrt{\frac{5N^2}{F(m+1, \beta)} + 4} \right) - \log\left(\frac{2\beta}{t_p}\right)}{2 \log \beta} \right\rfloor\\
	&=\left\lfloor m+1+\frac{\log (\beta/t_p)}{2\log \beta}       \right\rfloor\leq m+1.
\end{align*}
Similar arguments for $k_2^{+}$ and $k_2^{-}$ yield that 
\begin{align*}
	m&\leq k_2^+ \leq m+1,\\	
	m+1&\leq k_2^- \leq m+1.	
\end{align*}
In summary, since $s_p=k_1^+ + k_1^- + k_2^+ + k_2^- + 2$, for $p\in A_{m1}$ we have
\[
4m+4\leq s_p \leq 4m + 6,
\]
completing the proof of (i).

For the other three cases, one verifies that 
\begin{align*}
	&G(m+1, u_{m3})=F(m+1, u_{m3})<p\leq F(m+1, 1)<G(m+1,\beta),\ \textrm{if}\ p\in A_{m2};\\
	&G(m+1, 1)<F(m+1, 1)<p\leq F(m+1, u_{m1})=G(m, u_{m1}),\ \textrm{if}\ p\in A_{m3};\\
	&G(m, u_{m1})=F(m+1, u_{m1})<p\leq F(m, \beta)<G(m,1),\ \textrm{if}\ p\in A_{m4}.
\end{align*}
Then by the same method as that in the proof of (i), we get the desired results.\zb\\

We see from Theorem \ref{numofsolutions} that for a fixed $N$, the number of solutions to Equation \ref{eq:main} within the region $E_c(N)$ depends on the relative location of $p$ in $(0, N^2]$ rather than the value of $p$. One can verify that the estimate is sharp by computing with a concrete $N$.

\subsection{The estimate of $\mathbb{P}_{c_0}^{\mathcal{S}}(N^2)$ for $c_0=5/49$} As we can see from Theorem \ref{numofsolutions}, as $c$ decreases, the number of solutions to Equation \eqref{eq:main} within $E_c(N)$ increases. The point $c_0 = \frac{5}{49}$ is chosen so that $c_0N^2$ is the left endpoint of $A_{13}$, namely $c_0N^2 = F(2,1)$. Then the number of solutions to \eqref{eq:main} within $E_{c_0}(N)$ is more than $6$ for any $p$. We have the following estimate of $\mathbb{P}_{c_0}^{\mathcal{S}}(N^2)$.

\begin{thm}\label{densityforc0}
For \(c_0 = \frac{5}{49}\), it holds that
\begin{align*}
	\frac{3.6196}{\log N^{2}}+o\!\left(\frac{1}{\log N^{2}}\right)\leq\mathbb{P}_{c_0}^{\mathcal{S}}(N^2) \leq\frac{4.4536}{\log N^{2}}+o\!\left(\frac{1}{\log N^{2}}\right).
\end{align*}
\end{thm}	

\pf Since \(c_0 N^2 = F(2,1)\), another application of Theorem \ref{numofsolutions} gives
\begin{align*}
	\#\Omega_{c_0}^Y(N) \geq{}& \frac{7}{2}\big[\pi(F(2,1))-\pi(F(2,u_{13}))\big]+ \frac{8}{2}\big[\pi(F(2,u_{13}))-\pi(F(2,\beta))\big]\\	
	&+\frac{1}{2} \sum_{m=2}^{\left\lfloor\frac{\log N}{2\log\beta}\right\rfloor}
	\Big[
	(4m+1)\big(\pi(F(m+1,\beta^{-1})) - \pi(F(m+1,u_{m1}))\big) \\
	&+ (4m+2)\big(\pi(F(m+1,u_{m1})) - \pi(F(m+1,u_{m2}))\big) \\
	&+ (4m+3)\big(\pi(F(m+1,u_{m2})) - \pi(F(m+1,u_{m3}))\big) \\
	&+ (4m+4)\big(\pi(F(m+1,u_{m3})) - \pi(F(m+1,\beta))\big)
	\Big] + o\!\left(\frac{N^{2}}{\log N^2}\right)\\
	={}& \frac{1}{2\log N^2}\big[7 F(2,1)+ F(2,u_{13}) + \sum_{m=2}^{\left\lfloor\frac{\log N}{2\log\beta}\right\rfloor}
	\big[ F(m+1,\beta^{-1})+F(m+1,u_{m1})\\
	&+ F(m+1,u_{m2})+F(m+1,u_{m3}) \big]\big]+ o\!\left(\frac{N^{2}}{\log N^2}\right)\\
	\geq{}&0.4428\frac{N^2}{\log N^2}+o\!\left(\frac{N^{2}}{\log N^2}\right).
\end{align*}
Note that if $m>\left\lfloor\frac{\log N}{2\log\beta}\right\rfloor$, then $F(m+1, \beta^{-1})<2$, and there is no prime less than $2$. In view of Proposition \ref{XN}, we have 
\begin{align*}
	\mathbb{P}_{c_0}^{\mathcal{S}}(N^2) &= \frac{\#\Omega_{c_0}^Y(N)}{\# \Omega_{c_0}^X(N)} \geq\frac{3.6196}{\log N^2}+o\!\left(\frac{1}{\log N^2}\right).
\end{align*}

For the upper bound, using Theorem \ref{numofsolutions} again, we have 
\begin{align*}
	\#\Omega_{c_0}^Y(N)\leq{}&\frac{9}{2}\big[\pi(F(2,1))-\pi(F(2,u_{13}))\big]+ \frac{10}{2}\big[\pi(F(2,u_{13}))-\pi(F(2,\beta))\big]\\
	&+\frac{1}{2} \sum_{m=2}^{\left\lfloor\frac{\log N}{2\log\beta}\right\rfloor}
	\Big[
	(4m+3)\big(\pi(F(m+1,\beta^{-1})) - \pi(F(m+1,u_{m1}))\big) \\
	&+ (4m+4)\big(\pi(F(m+1,u_{m1})) - \pi(F(m+1,u_{m2}))\big) \\
	&+ (4m+5)\big(\pi(F(m+1,u_{m2})) - \pi(F(m+1,u_{m3}))\big) \\
	&+ (4m+6)\big(\pi(F(m+1,u_{m3})) - \pi(F(m+1,\beta))\big)
	\Big]+ o\!\left(\frac{N^{2}}{\log N^2}\right)\\
	={}& \frac{1}{2\log N^2}\big[9 F(2,1)+ F(2,u_{13}) + \sum_{m=2}^{\left\lfloor\frac{\log N}{2\log\beta}\right\rfloor}
	\big[ F(m+1,\beta^{-1})+F(m+1,u_{m1})\\
	&+ F(m+1,u_{m2})+F(m+1,u_{m3}) \big]\big]+ o\!\left(\frac{N^{2}}{\log N^2}\right)\\
	\leq{}&0.5449\frac{N^2}{\log N^2}+o\!\left(\frac{N^{2}}{\log N^2}\right).	
\end{align*}
It follows that 
\begin{align*}
\mathbb{P}_{c_0}^{\mathcal{S}}(N^2) =& \frac{\#\Omega_{c_0}^Y(N)}{\# \Omega_{c_0}^X(N)} \leq\frac{4.4536}{\log N^2}+o\!\left(\frac{1}{\log N^2}\right),
\end{align*}
and the proof is complete. \zb\\

Combining Theorem \ref{generalcase} and Theorem \ref{densityforc0}, we finally obtain the main result.

\begin{thm}\label{main1}
	For $c_0=5/49$, the following hold for the eigensurface ${\mathcal S}$ of the group ${\bf B}_3$:
	
	\emph{a)} $\liminf\limits_{N\to\infty} \frac{\mathbb{P}_{c_0}^\mathcal{S}(N)}{\mathbb{P}_{c_0}(N)}\geq  1.0037$;
	
	\emph{b)} $\limsup\limits_{N\to\infty} \frac{\mathbb{P}_{c_0}^\mathcal{S}(N)}{\mathbb{P}_{c_0}(N)}\leq  1.2350$.
\end{thm}

We note that the gap between the $\liminf$ and $\limsup$ of 
$\frac{\mathbb{P}_{c_0}^{\mathcal{S}}(N)}{\mathbb{P}_{c_0}(N)}$ 
depends on the estimation method employed. 
It would be of interest to investigate how to narrow this gap, 
and in particular, whether the limit itself exists. 
A central ingredient is the estimation, for each fixed positive integer $k$, 
of the proportion of positive primes $p<N$ such that Equation~\eqref{eq:main} 
has exactly $k$ solutions in the given region. 
To the best of our knowledge, this problem remains open.\\

\section{The limit behavior of $\mathbb{P}_{c}^{\mathcal{S}}(N^2)$}

This section examines the limit behavior of $\mathbb{P}_{c}^{\mathcal{S}}(N^2)$ as $c$ tends to $0$. For any fixed $0<c\leq 1$, as shown in Theorem \ref{main1}, it is possible to give an estimate of the limit superior and limit inferior of the ratio $\frac{\mathbb{P}_{c}^\mathcal{S}(N)}{\mathbb{P}_{c}(N)}$. It is thus meaningful to examine the limit behavior of the ratio $\frac{\mathbb{P}_{c}^\mathcal{S}(N)}{\mathbb{P}_{c}(N)}$ as $c$ tends to $0$. To our surprise, the following result holds.

\begin{thm}\label{main2}
The following limit exists:
\begin{align*}
	\lim_{c\to 0}\liminf_{N\to \infty} \frac{\mathbb{P}_{c}^\mathcal{S}(N)}{\mathbb{P}_{c}(N)}=\lim_{c\to 0}\limsup_{N\to \infty} \frac{\mathbb{P}_{c}^\mathcal{S}(N)}{\mathbb{P}_{c}(N)}=\frac{\sqrt{5}\,\pi^2}{36\,\zeta(3)\,\log\beta}\approx 1.0598.
\end{align*}
\end{thm}

\pf Recall from Section \ref{sec:solution_count} that the interval $(0,N^2]$ was first partitioned by $F(m,\beta)$, $m\geq 0$ and on each interval 
\begin{align*}
	A_m = (F(m+1, \beta), F(m, \beta)]=(F(m+1, \beta), F(m+1, \beta^{-1})],
\end{align*}
$F(m+1, u)$ is a monotonically decreasing function of $u$ on \([\beta^{-1}, \beta)\). Hence for any fixed $0<c\leq 1$, there is a unique positive integer $M$ and $u(c)\in [\beta^{-1}, \beta)$ such that $F(M+1, u(c))=cN^2$, where $M$ can be explicitly expressed as 
\begin{align}\label{M}
	M = \Biggl\lfloor 
	\frac{\log(\sqrt{5/c}+\sqrt{5/c-4})-\log(2\beta)}{2\log \beta}
	\Biggr\rfloor,
\end{align}
and for the number $u(c)$, only one of the following four cases will occur:
\begin{align*}\label{lowerb}
&\mathrm{(i)}\  u_{M3} \leq u(c) < \beta, \qquad\ \  \mathrm{(ii)}\ 1 \leq u(c) < u_{M3},\\
&\mathrm{(iii)}\ u_{M1} \leq u(c) < 1, \qquad \mathrm{(iv)}\ \beta^{-1} \leq u(c) < u_{M1}.
\end{align*}

Since the four cases are analogous, we only treat (i). In this case,
\begin{align*}
	\#\Omega_c^Y(N) \geq{}&\frac{1}{2}(4M+4)\big[\pi(F(M+1,u(c)))-\pi(F(M+1,\beta))\big]\\
	&+\frac{1}{2}\sum_{m=M+1}^{\left\lfloor\frac{\log N}{2\log\beta}\right\rfloor}  \Big[
	(4m+1)\big(\pi(F(m+1,\beta^{-1})) - \pi(F(m+1,u_{m1}))\big) \\
	&+ (4m+2)\big(\pi(F(m+1,u_{m1})) - \pi(F(m+1,u_{m2}))\big) \\
	&+ (4m+3)\big(\pi(F(m+1,u_{m2})) - \pi(F(m+1,u_{m3}))\big) \\
	&+ (4m+4)\big(\pi(F(m+1,u_{m3})) - \pi(F(m+1,\beta))\big)
	\Big]\\
	={}&\frac{1}{2\log N^2}\big[(4M+4)F(M+1,u(c))+             \sum_{m=M+1}^{\left\lfloor\frac{\log N}{2\log\beta}\right\rfloor}
	\big[ F(m+1,\beta^{-1} \\
	& +F(m+1,u_{m1}) + F(m+1,u_{m2})+F(m+1,u_{m3}) \big]\big]+o\!\left(\frac{N^{2}}{\log N^{2}}\right).
\end{align*}
On the other hand, 
\begin{align*}
\#\Omega_c^Y(N) 
\leq{}&\frac{1}{2\log N^2}\big[(4M+6)F(M+1,u(c))+             \sum_{m=M+1}^{\left\lfloor\frac{\log N}{2\log\beta}\right\rfloor}
\big[ F(m+1,\beta^{-1})\\
& +F(m+1,u_{m1})+ F(m+1,u_{m2})+F(m+1,u_{m3}) \big]\big]+o\!\left(\frac{N^{2}}{\log N^{2}}\right).
\end{align*}
For fixed $0<c<1$, we know from \eqref{M} that 
\begin{align*}
\frac{-\log c}{4\log \beta}-1 \leq M \leq  \frac{-\log c}{4\log \beta}+
	\frac{\log(\sqrt{5}/\beta)}{2\log \beta}.
\end{align*}
It follows that for fixed small $c$, when $N$ is large enough, the leading terms in $\#\Omega_c^Y(N)$ satisfy
\begin{align*}
	\frac{1}{2\log N^2}(4M+4)F(M+1,u(c))\sim\frac{-c\log c+O(c)}{2\log \beta}\cdot\frac{N^2}{\log N^2},\\
	\frac{1}{2\log N^2}(4M+6)F(M+1,u(c))\sim\frac{-c\log c+O(c)}{2\log \beta}\cdot\frac{N^2}{\log N^2}.
\end{align*}
The summation
\begin{align*}
   &\sum_{m=M+1}^{\left\lfloor\frac{\log N}{2\log\beta}\right\rfloor}
	\big[ F(m+1,\beta^{-1})+F(m+1,u_{m1}) + F(m+1,u_{m2})+F(m+1,u_{m3}) \big]\big]\\
	\leq{}&\sum_{m=M+1}^{\infty} 4F(m+1,\beta^{-1})\leq \frac{4}{\beta^4-1}cN^2=O(c)N^2.
\end{align*}
Hence,
\[
\#\Omega_c^Y(N)=\frac{-c\log c+O(c)}{2\log\beta}\cdot\frac{N^{2}}{\log N^{2}}+o\!\left(\frac{N^{2}}{\log N^2}\right).
\]
In view of Proposition \ref{XN}, we have 
\[
\mathbb{P}_c^{\mathcal{S}}(N^2)=\frac{\#\Omega_c^Y(N)}{\#\Omega_c^X(N)}= \frac{(-c\log c+O(c))/(2\log\beta)}{\frac{6}{\pi^2}\left(1-\int^1_{\frac{1-\sqrt{5-4c}}{2}}
		\left(1-\frac{\sqrt{5n^{2}+4c}-n}{2}\right)dn\right)}\frac{1}{\log N^2}+o\!\left(\frac{1}{\log N^2}\right).
\]
Therefore, for sufficiently small $c$,
\begin{align*} 
\lim_{N\to \infty} \frac{\mathbb{P}_{c}^\mathcal{S}(N)}{\mathbb{P}_{c}(N)}=\frac{\pi^2(-c\log c+O(c))}{36\zeta(3)\log\beta\left(1-\int^1_{\frac{1-\sqrt{5-4c}}{2}}
	\left(1-\frac{\sqrt{5n^{2}+4c}-n}{2}\right)dn\right) }.
\end{align*}
Letting \(c\to 0\) gives
\[
\lim_{c\to 0}\liminf_{N\to \infty} \frac{\mathbb{P}_{c}^\mathcal{S}(N)}{\mathbb{P}_{c}(N)}=\lim_{c\to 0}\limsup_{N\to \infty} \frac{\mathbb{P}_{c}^\mathcal{S}(N)}{\mathbb{P}_{c}(N)}=\frac{\sqrt{5}\pi^2}{36\zeta(3)\log\beta}\approx 1.0598.
\]
This completes the proof of Theorem \ref{main2}.\zb\\

\section{Concluding remark}

This paper shows that the eigensurface ${\mathcal S}$ acts as a ``prime-enhancing" sieve, selecting integer points that are more likely to contain a prime coordinate. This suggests that the eigensurface of a group representation carries more than just algebraic information--it may encode subtle ``arithmetic signature" that affects the distribution of primes. The work here naturally points to three promising directions for future research:

{\bf 1}. The methodology developed here for the braid group ${\bf B}_3$ and its reduced Burau representation may be extended to other finite groups and representations. Similar study can also be carried out for classical simple Lie algebras, whose characteristic polynomials have been studied intensely in recent years \cite{GLW,Ya2}. The complexity of the Diophantine equations involved will increase, and this will help us gain deeper number-theoretic insights on representation theory. 

{\bf 2}.  Beyond the braid group ${\bf B_3}$, this ``prime-enhancing" phenomenon has also been observed on the symmetry group ${\bf S}_3$ \cite{GHY} and the free group ${\bf F}_2$. However, the underlying mechanism remains mysterious. Therefore, a central open problem is to identify the exact structural properties of the groups and their representations that lead to this phenomenon. 

{\bf 3}. The results of this paper may have implications for the study of arithmetic dynamics, as shown in the recursive structure of the solutions to the equation $m^2+mn-n^2=p$. The eigensurfaces of group representations serve as a bridge connecting group dynamics and arithmetic dynamics. Perhaps this connection holds the key to the mystery.
 
In conclusion, this paper is one of a series of works in the nascent field of ``arithmetic representation theory," which sits at the intersection of representation theory, arithmetic geometry, analytic number theory, and spectral theory. It contributes both meaningful results and a broader research program that promises to enrich all of its constituent disciplines. \\

\noindent\textbf{Conflict of Interest.} The authors declare that they have no conflict of interest.\\

\noindent\textbf{Acknowledgments.} This work is partially supported by the National Natural Science Foundation of China under Grant No. 12271090 and the Jiangsu Provincial Scientific Research Center of Applied Mathematics under Grant No. BK20233002. The fourth author is supported partially by the Simons Foundation Grant No. MP-TSM-00002315.

\end{document}